# Least Squares Fitting of Low-Level Gamma-ray Spectra with B-Spline Basis Functions

M.H. Zhu*, L.G. Liu, Z. You, A.A. Xu

Space Exploration Technology Laboratory,
Macao University of Science and Technology, Avenue Wai Long, Taipa, Macao

**Abstract:** In this paper, new methods for smoothing gamma-ray spectra measured by NaI detector are derived. Least squares fitting method with B-spline basis functions is used to reduce the influence of statistical fluctuations. The derived procedures are simple and automatic. The results show that this method is better than traditional method with a more complete reduction of statistical fluctuation.

## 1. Introduction

The objective of analysis of gamma-ray spectra is to extract the useful information from the measured data while minimizing the influence of statistical fluctuations that can raise calculation error on both nuclide recognition and quantitative analysis. A number of mathematical procedures have been developed for processing detected data and extracting useful information, for example, such as Polynomial Fitting[1~4], Fourier transformation[5~9], and Convolution operation[10~14] which is the most well known and important method for data smoothing[15].

However, applications of these methods are not always efficacious, especially in the case of intensive noised spectra as Fig.1. The polynomial fitting method applied to remove the statistical fluctuation always comes along with some problems that can be also found in Convolution operation method, for example, spectral distortion, weak peaks easily lost and false peaks generated that then can rise calculation error in background determination, peak searching , fitting, etc. The Fast Fourier Transformation (FFT) method has little application in this area because of little energy information reserved in the transformed data[16], and also too much computer time consumed. Some math filters, for example, such as the Gaussian distribution function, the Lorentzian distribution function are used in Convolution operation method to remove statistical fluctuations. Better result will be obtained while the value of math filter's full width at half maximum (FWHM) is approximated to the average FWHM(in channels) of the peaks in the spectra to be smoothed. However, it is hard to calculate the average FWHM in the fluctuated spectra. Rough value used in the processing always gives uncompleted elimination result and cannot be fit for different spectrum. Sometimes, a simple approximated method that convolute the spectra with a vector in the form of ***n***-points weight means[17] , for example, such as $[1, 2, 1]^T/4$ for three-point weighted means , $[-3, 12, 17, 12 , -3]^T/35$ for five- point weighted means , is used to replace the math filter functions to reduce the statistical fluctuations . These replacements also have some disadvantages: smoothing is too intensive to flat the peaks in some portions while in significant residual fluctuation smoothing is too mild[16].

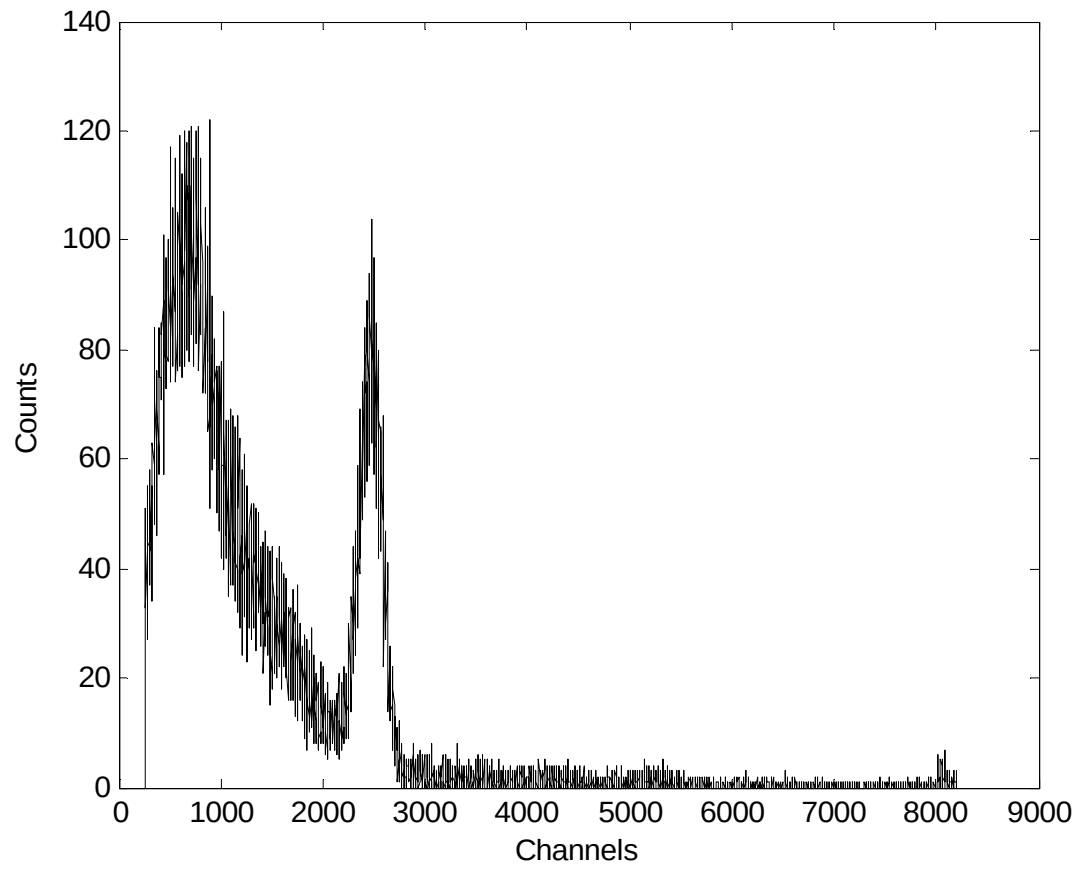

Fig.1. gamma-ray spectrum of $^{137}$Cs (NaI)

*Tel.: +00853 8972024;
E-mail address: peter_zu@163.com

In this paper, least squares method with B-spline basis functions is introduced to fit the gamma-ray spectra without re-processing. The description of this method and its application is divided into four sections. In section 2, B-spline basis functions and least squares are described briefly. In section 3, gamma-ray spectra measured in the laboratory are used to test this method. Conclusions will be given in the last section.

## 2. Description of the method

The B-spline basis functions is defined as a piecewise function which is non-zero only over four adjacent intervals between knots.

$$\varphi(x) = \begin{cases} 0 & |x| \geq 2 \\ \frac{1}{2}|x|^3 - x^2 + \frac{2}{3} & |x| \leq 1 \\ -\frac{1}{6}|x|^3 + x^2 - 2|x| + \frac{4}{3} & 1 < |x| < 2 \end{cases}$$

Fig.2 typifies the shape of this function. Because it is a spline and therefore has the appropriate continuity at given knots, its value and first and second derivatives are zero at these knots[18].

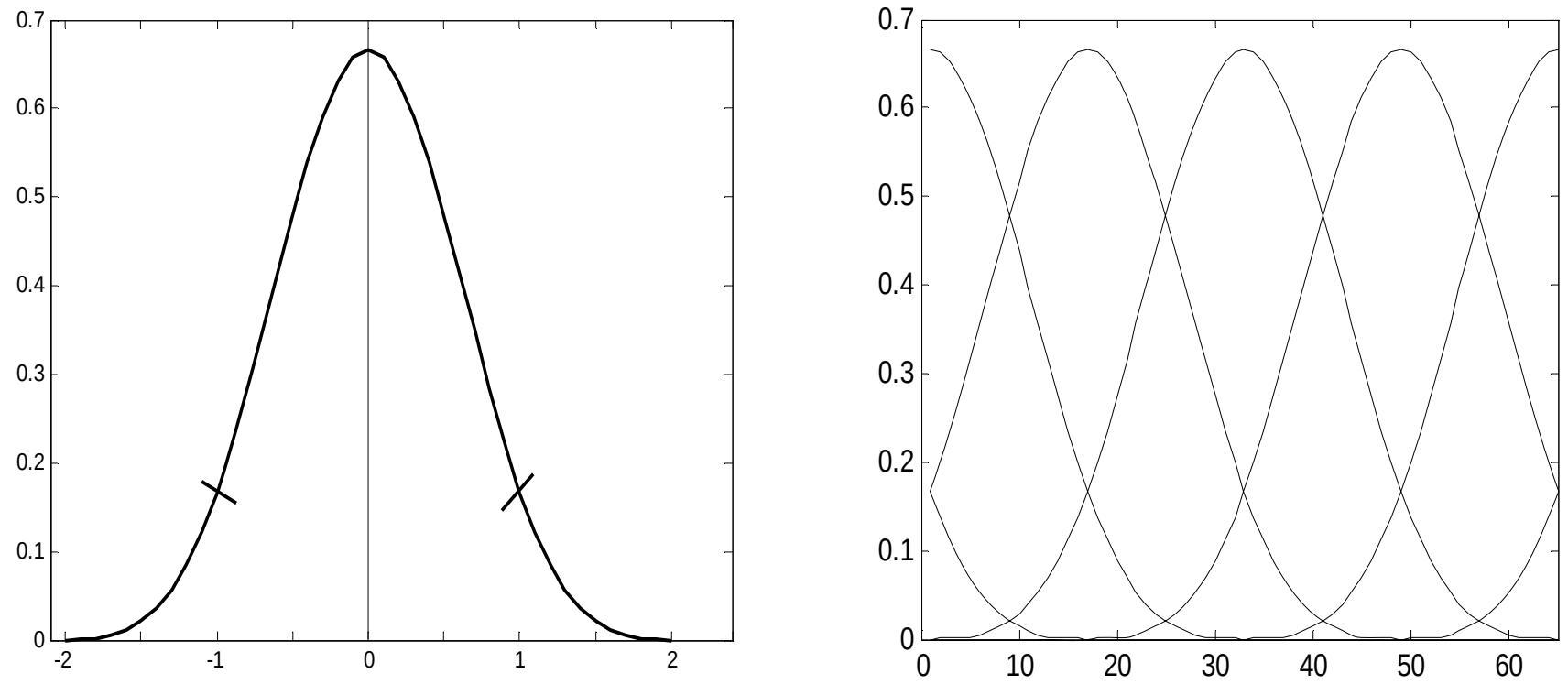

(a) (b)

Fig.2. (a) B-splilne basis functions shape
(b) 5-fundamental splines ( $c_i$=1 )

To define the full set of B-splines which are required for our purpose in the range of interest channels $a \leq x \leq b$, where $x$ is channel number. Firstly, it is necessary to select five channels $x_0, x_1, x_2, x_3, x_4$ with equal interval that satisfy $a \leq x_0 < x_1 < x_2 < x_3 < x_4 \leq b$. With this set of channels, define, as above, the fundamental splines $\varphi_j(x)$, $j = 1, 2, \ldots, 5$. Than the general B-splines with channels $x_0, x_1, x_2, x_3, x_4$ has the unique representation in the range $a \le x \le b$ of the form

$$s(x) = \sum_{j=1}^{5} c_j \varphi_j(x).$$

Fig.2(b) shows the shape of 5-fundamental splines while $c_j = 1$. The method of least squares assumes that the best-fit $s^*(x)$ curve with these selected channels has the minimal sum of the least square error from the measured gamma-ray spectra $N(x)$.

$$\left\|\delta_i\right\|^2 = \sum_{x=1}^{m} [S_i^*(x) - N(x)]^2 = \min \sum_{x=1}^{m} [S_i(x) - N(x)]^2 ,$$

$$S_i(x) = C_{i,0}\varphi_0(x) + C_{i,1}\varphi_1(x) + \cdots + C_{i,n}\varphi_n(x) \qquad n < m ,$$

Where ***m*** is the number of the channels of the spectra in the interval [a, b], and ***n*** is the number of the selected channels. Then, new channel is selected as knots in the middle of the interval and the process is repeated until the interval is equal to 1. Best-fit function $s_i^*(x)$ with selected channels will be obtained at each time. The corresponding best-fit curves also have minimal sum of the least square error each other

$$\left\|\varepsilon_i\right\|^2 = \sum_{all\,x} [S_{i+1}^* - S_i^*]^2 .$$

The result curve $s_i^*(x)$ will be obtained while $\left\|\varepsilon_i\right\|^2$ gets the global minimum.

## 3. Test

As a test of this method, the gamma-ray spectrum as shown in Fig.1 is measured from $^{137}Cs$ source by NaI detector for more than 2 hours. The spectra, processed by B-spline basis functions fitting and traditional 3-points weighted averaged smoothing (iterates 5000), are respectively given in Fig.3(a) and(b). It can be seen from the figures that, the B-spline basis functions method can eliminate statistical fluctuation thoroughly with little data distortion and the result best-fit curve is smooth compared with the traditional 3 or 5-points weighted means smoothing method. This is of important significance in peak searching and distinguishing the overlapping peaks.

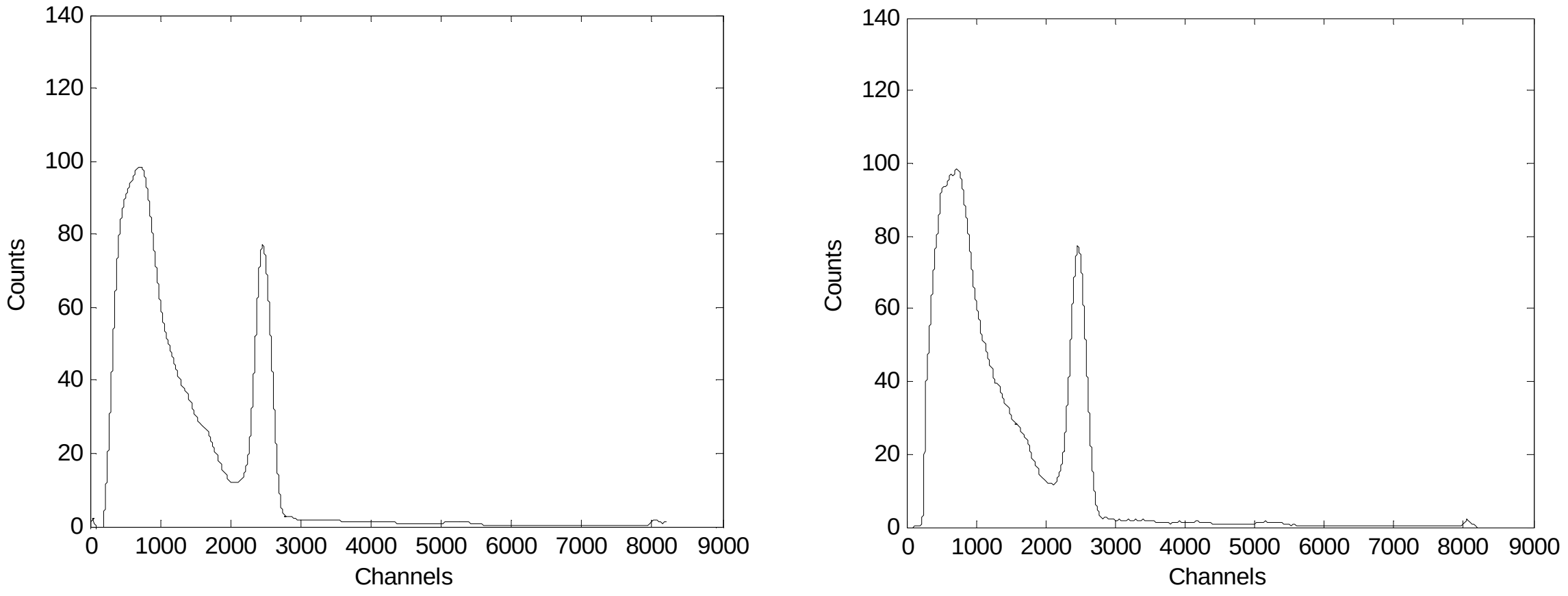

(a) (b)

Fig.3. (a) Smoothed data of $^{137}Cs$ (B-Spline basis functions)
(b) Smoothed data of $^{137}Cs$ (3-points weighed average – 5000 Iterations)

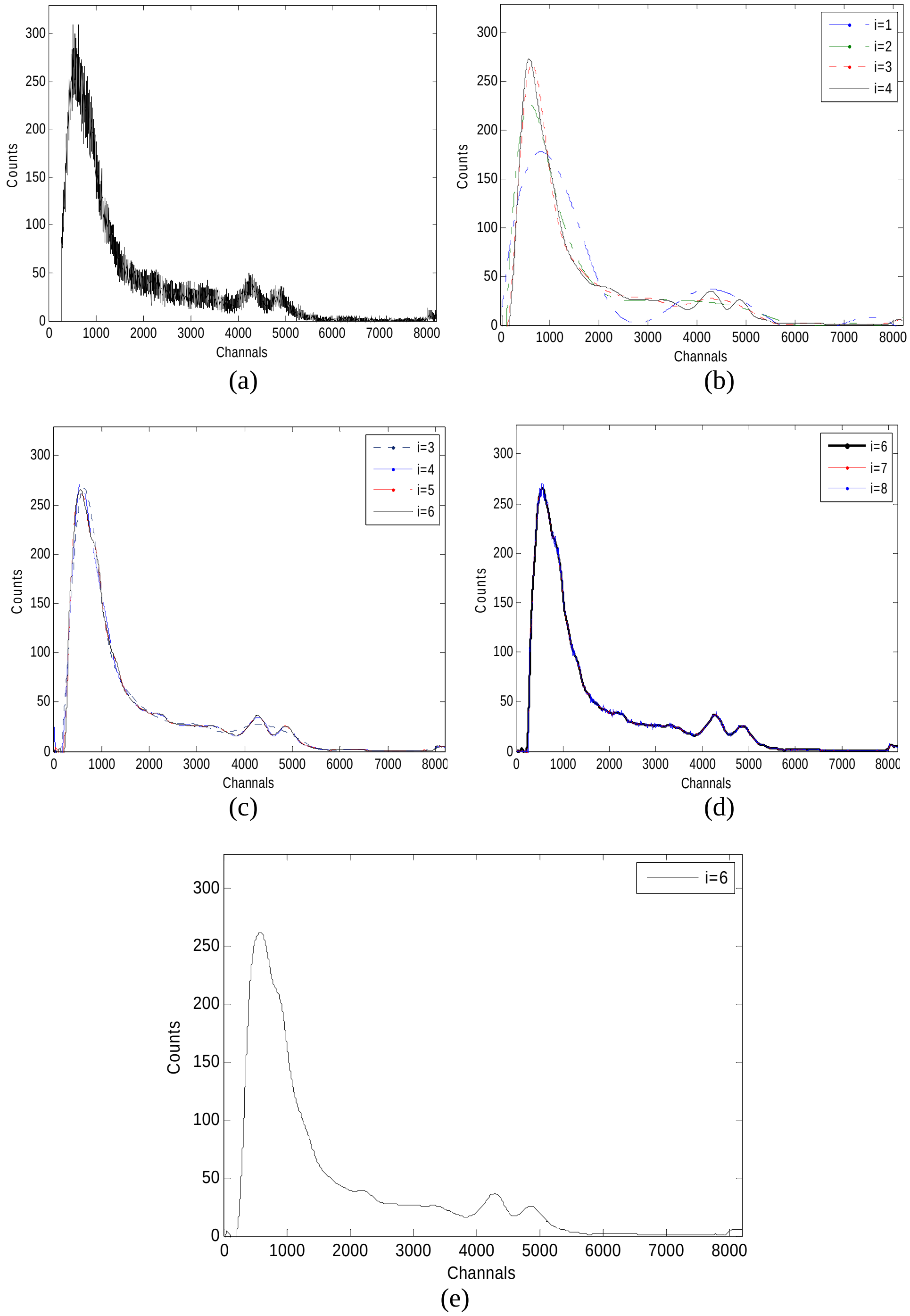

Fig.4. (a) Gamma-ray spectrum of $^{60}$Co (NaI)

(b) $S_i^*$ ($\boldsymbol{i}$ = 1 ~ 4) (c) $S_i^*$ ($\boldsymbol{i}$ = 3 ~ 6) (d) $S_i^*$ ($\boldsymbol{i}$ = 6 ~ 8) (e) The result fit curve while $\boldsymbol{i}$ = 6.

Gamma-ray spectrum of $^{60}$Co and the mixed radiate source are also measured by NaI detector to test this method, as shown in Fig.4 and Fig.5, respectively. In Fig.4, best-fitting curve functions $s_i^*(x)$ are also shown with subscript ***i*** increasing from 1 to 8. The results of $\|\varepsilon_i\|^2$, calculated, as above, are listed in Table 1. It can be seen from the table that $\|\varepsilon_i\|^2$ gets some minimum values when subscript ***i*** increases. The best-fitting curve is selected as result from $s_i^*(x)$ while $\|\varepsilon_i\|^2$ gets the global minimum value shown in bold type.

Table.1 Least square error of corresponding best-fit curves ($10^5$)

| **Source** | **Channels** | $\|\varepsilon_1\|^2$ | $\|\varepsilon_2\|^2$ | $\|\varepsilon_3\|^2$ | $\|\varepsilon_4\|^2$ | $\|\varepsilon_5\|^2$ | $\|\varepsilon_6\|^2$ | $\|\varepsilon_7\|^2$ | $\|\varepsilon_8\|^2$ |
|---|---|---|---|---|---|---|---|---|---|
| $^{137}$Cs | 8192 | 4.3814 | 0.4671 | 0.1991 | 0.0583 | 0.0247 | ***0.0036*** | 0.0037 | 0.0027 |
| $^{60}$Co | 8192 | 3.6636 | 0.4345 | 0.1802 | 0.0263 | 0.0164 | ***0.0035*** | 0.0036 | 0.0026 |
| Mixture | 8192 | 8.7579 | *3.7614* | 6.0633 | 5.6871 | 0.9147 | ***0.1464*** | 0.1969 | 0.3913 |

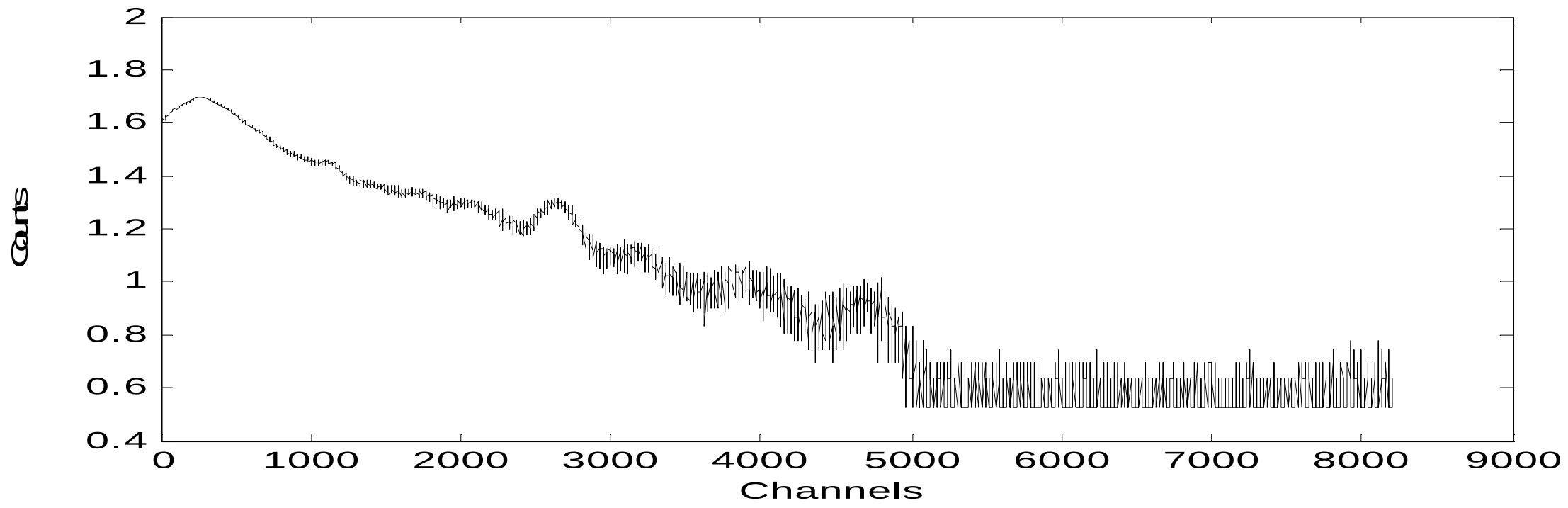

(a)

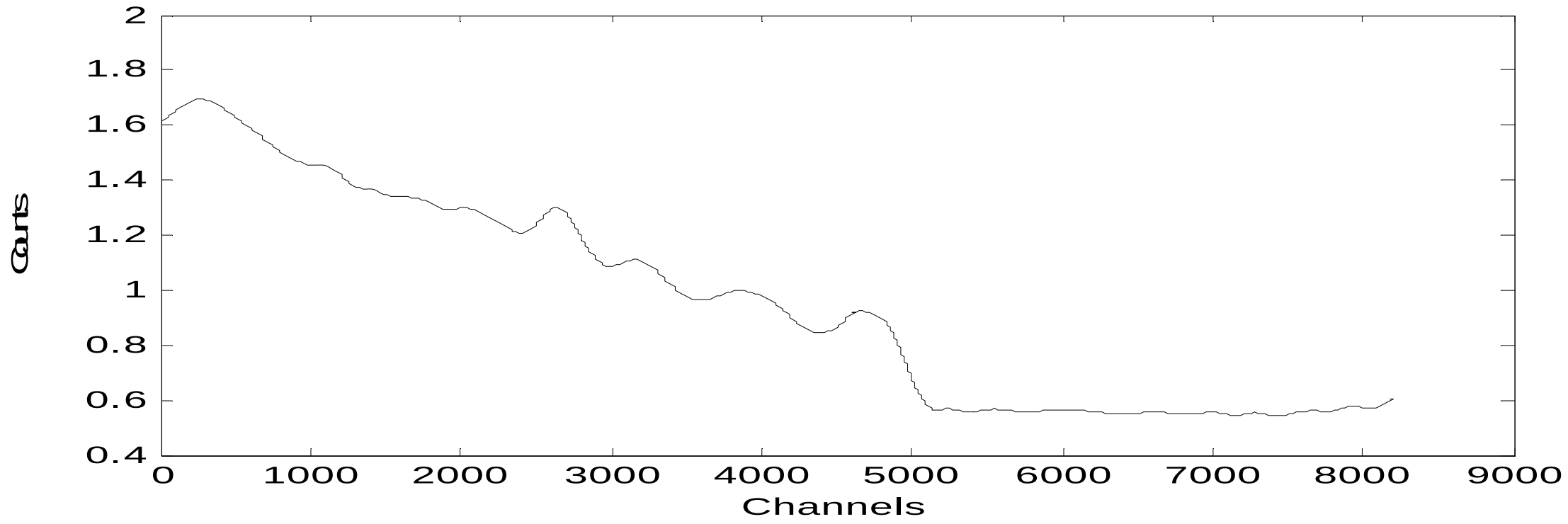

(b)

Fig.5 (a) Mixed radiate source (NaI)
(b) Smoothed data of mixed radiate source using B-Spline basis functions

# 4. Conclusion

In the paper, a least squares fitting method with B-splines basis functions has been developed for eliminating the statistical fluctuation of gamma-ray spectra measured by NaI detector. As one can see from the calculations as above, the method is simple and does not need any initial input value. Different tests show that this method can remove statistical fluctuation completely with little data distortion and the result best-fit curve is smooth which it is important in the background elimination, peak searching and distinguishing the overlapping peaks in the analysis of gamma-ray spectra.

The authors would like to thank Professor DongXu Qi for his useful advice and Jian Li for useful discussion. Financial support of the Science and Technology of Development Fund of Macao is gratefully acknowledged.